\documentclass[a4paper]{article}

\usepackage{amsmath}
\usepackage{amsfonts}
\usepackage{amssymb}
\usepackage{latexsym}
\usepackage{epsfig} 
\usepackage{psfrag}

\scrollmode

\newcommand{\begriff}[1]{\textbf{#1}}

\renewcommand{\d}{\partial}
\renewcommand{\t}{\mathcal T}
\newcommand{\C}{\mathcal C}
\newcommand{\p}{\mathbb P^3}
\newcommand{\z}{\mathcal Z}

\newcommand{\co}{\colon\thinspace}    

\newtheorem{theorem}{Theorem}

\newtheorem{lemma}{Lemma}

\newtheorem{definition}{Definition}

\def\qed{\relax\ifmmode\expandafter\endproofmath\else
  \unskip\nobreak\hfil\penalty50\hskip.75em\hbox{}\nobreak\hfil\bull
  {\parfillskip=0pt \finalhyphendemerits=0 \bigbreak}\fi}
\def\endproofmath$${\eqno\bull$$\bigbreak}
\def\bull{\vbox{\hrule\hbox{\vrule\kern3pt\vbox{\kern6pt}\kern3pt\vrule}\hrule}}

\newcounter{constr}

\begin{document}

\title{Complexity of triangulations of the projective space}  

\author{Simon King\thanks{ Department of Mathematics,
    Darmstadt University of Technology, Schlossgartenstr.~7, 64289~
    Darmstadt, Germany.  E-mail: \texttt{king@mathematik.tu-darmstadt.de}}}

\date{\today}
\maketitle

\begin{abstract}
  It is known that any two triangulations of a compact $3$--manifold
  are related by finite sequences of certain local transformations. We
  prove here an upper bound for the length of a shortest
  transformation sequence relating any two triangulations of the
  $3$--dimensional projective space, in terms of the number of
  tetrahedra.

\medskip
\noindent
\textbf{MSC--class:} 57Q25; 57Q15.
\end{abstract}

\section{Introduction}
\label{sec:intro}

By the \lq\lq Hauptvermutung\rq\rq, that was proven by
Moise~\cite{moise}, any two triangulations $\t_1$ and $\t_2$ of a
compact $3$--manifold $M$ have a common subdivision. This allows to show
that $\t_1$ and $\t_2$ are related by finite sequences of certain local
transformations of triangulations, e.g., stellar
subdivisions~\cite{alexander} or elementary shellings~\cite{pachner}.
These results do not provide explicit constructions of transformation
sequences and do not yield a recognition algorithm for $M$.
In this paper, we construct transformation sequences for
triangulations of the $3$--dimensional projective space $\p$.  We
consider the following local transformations, that generalise stellar
subdivisions.

\begin{definition}
  Let $\t$ and $\tilde\t$ be PL--triangulations of a closed PL--manifold, and 
  let  $e$ be an edge of
  $\t$ with $\d e=\{a,b\}$. Suppose that $\tilde\t$ is obtained from
  $\t$ by removing the open star of $e$ and identifying $a*\sigma$
  with $b*\sigma$ for any simplex $\sigma$ in the link of $e$.
  Then $\tilde\t$ is the result of the \begriff{edge contraction} of
  $\t$ along $e$, and $\t$ is the result an \begriff{edge expansion} of
  $\tilde\t$ along $e$. 
\end{definition}
\noindent
In general, there are edges of $\t$ along which a contraction is
impossible. This is the case, e.g., if the edge is part
of an edge path of length 3 that does not bound a 2--simplex of
$\t$. 
It is easy to see that any PL--triangulation admits only a finite number of
edge expansions. 

Let $d(\t_1,\t_2)$ be the length of a shortest sequence of edge
contractions and expansions relating two triangulations $\t_1$ and
$\t_2$ of a closed $3$--manifold $M$.
The aim of this paper is to provide an upper bound for $d(\t_1,\t_2)$ when
$M$ is homeomorphic to $\p$, as stated in
Theorem~\ref{thm:main} below.
The proof is partially based on our work on the $3$--sphere 
(see~\cite{king1}--\cite{king3}).
This paper is thought of as a first step towards a study of more general
$3$--dimensional manifolds, e.g., (atoroidal) Haken manifolds.
A generalisation to \emph{all} compact $3$--manifolds, which would solve the
algorithmic classification problem for compact $3$--manifolds, is out of
reach, as yet.

\begin{theorem}\label{thm:main}
  Any two triangulations of\/ $\p$ with at most $t$ tetrahedra are
  related by a sequence of less than $2^{27000\,t^2}$ edge
  contractions and expansions.
\end{theorem}
The constant factor in the exponent is certainly not optimal.
According to the examples in~\cite{king2}, concerning the minimal number
of edge expansions needed to transform a triangulation of $S^3$ into a
\emph{polytopal} triangulation, we believe that the bound in
Theorem~\ref{thm:main} can not be replaced by a subexponential bound.

We outline the proof of Theorem~\ref{thm:main}. It is based on Haken's
normal surface theory, Barnette's work~\cite{barnette} on irreducible
triangulations of the projective plane, and our techniques
in~\cite{king1} and~\cite{king2}.
Let $\t_1$ be a triangulation of $\p$.
By means of normal surface theory, we construct a certain projective
plane $P\subset \p$.  The complement of a regular neighbourhood of $P$
is a ball. This allows to apply techniques of~\cite{king1}
and~\cite{king2}, yielding a sequence of edge contractions and
expansions relating $\t_1$ with a triangulation $\t_P$ that depends on
the choice of $P$.
The next step is to simplify $\t_P$ by a certain series of edge
contractions.  In that way we transform $\t_1$ via $\t_P$ into one of
two \emph{standard triangulations} of $\p$, corresponding to the two
irreducible triangulations of the projective plane found by
Barnette~\cite{barnette}.
In the last step in the proof of Theorem~\ref{thm:main}, we relate the
two standard triangulations of $\p$ by an explicit sequence of edge
contractions and expansions. In conclusion, we construct
transformation sequences that relate any two triangulations of $\p$
via the two standard triangulations. The bound stated in
Theorem~\ref{thm:main} follows from a complexity analysis of the
construction.

The paper is organised as follows.  In Section~\ref{sec:normal}, we
briefly outline the results of normal surface theory used in our
proofs. Section~\ref{sec:contractions} recalls a lemma
from~\cite{king2} on the construction of sequences of edge contractions. 
We prove the existence of $P$ in
Section~\ref{sec:gibtP}.  Section~\ref{sec:intoTP} is devoted to the
transformation into $\t_P$, along the lines of~\cite{king1}--\cite{king3}.
Section~\ref{sec:standardtriang} is concerned with the transformation into
one of the two standard triangulations. A transformation sequence
relating the two standard triangulations is finally given in
Section~\ref{sec:proofend}.

\section{Prerequisites}
\label{sec:voraussetzungen}

In this section, we collect some results used in the proof
of Theorem~\ref{thm:main}. 
We denote the number of connected components of a topological space $X$
by $\#(X)$.  For a tame subset $Y\subset X$, we denote an open regular
neighbourhood of $Y$ in $X$ by $U(Y)$.

\subsection{Normal surfaces}
\label{sec:normal}

Let $M$ be a closed 3--manifold with a cellular decomposition $\z$, so
that the closure of any open cell of $\z$ is homeomorphic to a closed
ball. Its $k$--skeleton is denoted by $\z^k$. In our applications, $\z$
is a triangulation or is dual to a triangulation. 
A \begriff{normal isotopy} with respect to $\z$ is an ambient isotopy of
$M$ that fixes each cell of $\z$ set-wise.  
\begin{definition}
  Let $c$ be a closed $2$--cell of $\z$ and let $\gamma\subset c$ be a
  closed embedded arc with $\gamma\cap\d c=\d\gamma$, disjoint from
  the vertices of $c$. If $\gamma$ connects two different edges in the
  boundary of $c$ then $\gamma$ is a \begriff{normal arc}.  Otherwise
  it is a \begriff{return}.
\end{definition}

\begin{definition}\label{def:normal}
  Let $S\subset M$ be a closed embedded surface transversal to $\z$.  We
  call $S$ \begriff{normal} with respect to $\z$, if $S\cap\z^2$ is a
  union of normal arcs, $S\setminus\z^2$ is a disjoint union of discs,
  and the boundary of any connected component of $S\setminus\z^2$ meets
  any edge of $\z$ at most once.
\end{definition}

\begin{figure}[htbp]
  \begin{center}
    \leavevmode
    \epsfig{width=7cm,file=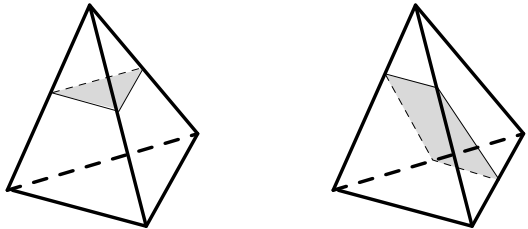}
    \caption{Normal triangles and squares}
    \label{fig:pieces}
  \end{center}
\end{figure}

When it is clear from the context, we do not specify with respect to
which cellular decomposition a surface is normal. In the rest of this
section, we focus on normal surfaces with respect to a triangulation
$\t$ of $M$. 
It is well known that normal surfaces in a triangulated $3$--manifold
are built from copies of so-called \emph{normal triangles} and
\emph{normal squares} (see Figure~\ref{fig:pieces}).

\begin{theorem}\label{thm:kneser}
  Let $\t$ be a triangulation of $M$ with $n$ tetrahedra.  Let
  $S\subset M$ be a normal surface comprising
  more than $10 n$ two-sided connected components. Then two connected
  components of $S$ are normally isotopic.\qed
\end{theorem}
This result is originally due to Kneser and proven, e.g., in Lemma~4
of~\cite{haken2}.  We use this result to show that certain iterative
constructions of normal surfaces stop after a finite number of steps.

Under a technical condition (see~\cite{hemion} for details), one can
define the \emph{sum} $S_1+S_2$ of two normal surfaces $S_1,S_2\subset
M$.  The sum is a normal surface and is determined up to normal
isotopy by the condition $(S_1+ S_2)\cap \t^1 = (S_1\cup S_2)\cap
\t^1$. The Euler characteristic is additive, $\chi(S_1+S_2) =
\chi(S_1)+\chi(S_2)$.
Haken has shown that with this notion of a sum, the set of normal
surfaces in $M$ with respect to $\t$ is isomorphic to a subgroupoid of
the semi group $G$ of integer points in a rational convex cone, the
so-called Haken cone.
The semi group $G$ is additively generated by a finite set of
elements, that can be constructed by means of integer programming.
The set of normal surfaces is finitely generated as well, by the
following result. For an embedded surface $S\subset M$ that is in
general position with respect to $\t$, denote $\|S\|= \#(S\cap \t^1)$.

\begin{theorem}\label{thm:fundamental}
  Let $N\subset M\setminus U(\t^0)$ be a sub--$3$--manifold whose
  boundary is a normal surface.  There is a system $F_1,\dots,
  F_q\subset N$ of normal surfaces such that $$\|F_i\| < \|\d N\|\cdot
  2^{18 n}$$ for $i=1,\dots,q$, and any normal surface $F\subset N$ can
  be expressed as a sum $F=\sum_{i=1}^q k_iF_i$ with non-negative
  integers $k_1,\dots,k_q$. \qed
\end{theorem}
The surfaces $F_1,\dots,F_q$ are called \begriff{fundamental surfaces
  in $N$}.
The preceding theorem is Theorem~3 of~\cite{king1}, that is formulated
in a slightly more general setting, namely for so-called $2$--normal
surfaces. One can prove slightly better bounds in our special
situation, but this concerns only the constant $18$.
Actually it follows from~\cite{hass} that the bound in
Theorem~\ref{thm:fundamental} can not be replaced by a subexponential
bound.
In the case $N=M\setminus U(\t^0)$, the existence and constructibility
of fundamental surfaces 
is classical in Haken theory~\cite{hemion}, and bounds for $\|F_i\|$
were obtained in~\cite{hasslagarias}.

\subsection{Edge contractions}
\label{sec:contractions}

We recall in this subsection a lemma yielding sequences of edge
contractions and expansions. 
Let $M$ be a closed $3$--manifold. A cellular decomposition of $M$ is
\begriff{simple}, if any vertex is adjacent to exactly four edges, and
any edge is adjacent to exactly three $2$--cells (counted with
multiplicity).  A cellular decomposition is \begriff{regular} if the
closure of any open $k$--cell is homeomorphic to a $k$--ball
($k=1,\dots,3$.

Let $\C_1$ be a simple regular cellular decomposition.  In
general, $\C_1$ is not dual to a triangulation, as multiple edges might
occur. However, the barycentric subdivision $\C_1'$ of $\C_1$ is a
triangulation of $M$.
By the next lemma, that is proven in~\cite{king2}, the deletion of an
appropriate $2$--cell of $\C_1$ gives rise to a sequence of 
contractions of $\C_1'$.
\begin{lemma}\label{lem:contr}
  Let $\C_1,\C_2$ be two simple regular cellular decompositions of
  $M$, so that $\C_1^2\setminus \C_2^2$ is an open $2$--cell of $\C_1$
  with $k$ vertices in its boundary.  Then $\C_2'$ is obtained from
  $\C_1'$ by a series of $4k+2$ contractions. \qed
\end{lemma}

\section{Transforming triangulations}
\label{sec:triangtrafo}

The four parts of this section form the proof of 
Theorem~\ref{thm:main}, as outlined in the introduction.
We fix the following notations.  Let $\C$ be a cellular decomposition
of $\p$ that is dual to a triangulation, and let $\t$ be its
barycentric subdivision. Let $n$ be the number of tetrahedra of $\t$.
We consider $\C^2$ as a subset of $\t^2$.

\subsection{Fundamental projective plane}
\label{sec:gibtP}

It is well known (see~\cite{hemion}) that there is a fundamental
projective plane with respect to $\t$. We prove here that we can choose
it so that additionally it is normal with respect to $\C$. This 
technical condition is needed to make the techniques of~\cite{king2} work 
(see next subsection).
\begin{lemma}\label{lem:gibtP}
  Among the fundamental surfaces with respect to $\t$ in $\p\setminus
  U(\t^0)$, there is a projective plane $P$ that is normal with
  respect to $\C$.
\end{lemma}
\begin{proof}
  Choose an embedded projective plane $P\subset \p$ so that the triple
  $$(\#(P\cap \C^1),\#(P\cap \t^1),\#(P\cap \t^2))$$
  is minimal in lexicographic order. 
  We first prove that $P$ is normal with respect to $\t$ and $\C$, by a
  modification of standard techniques (compare~\cite{hemion}).  Since
  $\#(P\cap \C^1)$ is minimal, $P\cap \C^2$ contains no returns.
  Removing a return in $P\cap \t^2$ does not increase $\#(P\cap \C^1)$,
  thus there is no return in $P\cap \t^2$ by minimality of $\#(P\cap
  \t^1)$.
  The minimality of $\#(P\cap \t^2)$ excludes circles in $(P\cap
  \t^2)\setminus \t^1$, since cutting-and-pasting along such a circle
  does not increase $\#(P\cap \C^1)$ and $\#(P\cap \t^1)$.
  
  Assume that there is a circle $\gamma$ in $(P\cap\C^2)\setminus
  \C^1$, contained in a $2$--cell $c$ of $\C$. By the preceding
  paragraph, we know that $\gamma$ is a union of normal
  arcs with respect to $\t$, and $\gamma$ bounds a disc in $c$
  containing a vertex of $\t$ (the barycenter of $c$).  Therefore
  $\#(\gamma\cap \t^1)\ge 6$.
  Let $D\subset c$ be the disc bounded by $\gamma$. Since $P$ is
  incompressible in $\p$, there is a disc $D'\subset P$ bounded by
  $\gamma$. We replace $D'$ by a parallel copy $D''$ of $D$ with $\d
  D''\cap D'=\emptyset$. We can choose $D''$ so that $\#(D''\cap
  \t^1)=1$, namely intersecting an edge of $\t$ that connects the
  barycenter of $c$ with the barycenter of a $3$--cell of $\C$. Hence,
  the cut-and-paste operation replacing $D'$ by $D''$ decreases
  $\#(P\cap \t^1)$ by at least $5$, without increasing $\#(P\cap
  \C^1)$. So by minimality of $\#(P\cap \t^1)$, there is no circle in
  $(P\cap\C^2)\setminus \C^1$.
  
  By the preceding paragraph, any boundary component of a connected
  component of $P\setminus \C^2$ meets $\C^1$. Thus the connected
  component of $P\setminus \C^2$ are discs, since otherwise one can
  decrease $\#(P\cap \C^1)$ by a cut-and-paste operation.
  If the boundary of a connected component $P\setminus \C^2$ intersects
  some edge of $\C$ at least twice, then there is a closed embedded disc
  $D\subset \p$ so that $D\cap \C^2 = \d D\cap \C^1$ is an arc contained
  in the interior of an edge of $\C$, $D\cap P\subset \d D$, and $\d D\subset
  P\cup \C^1$. Sliding $P$ across $D$ decreases $\#(P\cap\C^1)$.
  The corresponding operations, applied to connected components of
  $P\setminus \t^2$, decreases $\#(P\cap \t^1)$ without increasing
  $\#(P\cap \C^1)$. In conclusion, the minimality of $(\#(P\cap
  \C^1),\#(P\cap \t^1),\#(P\cap \t^2))$ implies that the boundary of
  any connected component of $P\setminus \C^2$ (resp. $P\setminus
  \t^2$) meets any edge of $\C$ (resp. of $\t$) at most once. Hence
  $P$ is normal both with respect to $\t$ and to $\C$.
  
  We represent $P$ as a sum $F_1+\dots + F_k$ of fundamental surfaces
  with respect to $\t$. We can assume that none of $F_1,\dots ,F_k$ is a
  $2$--sphere (see~\cite{hemion}). Since the Euler characteristic is
  additive under the addition of normal surfaces, one summand (say,
  $F_1$) has positive Euler characteristic, thus, is a projective plane.
  Since $\C^1\subset \t^1$, we have $\#(F_1\cap\C^1)\le \#(P\cap
  \C^1)$ and $\#(F_1\cap\t^1)\le \#(P\cap \t^1)$, and
  $\#(F_1\cap\t^2)\le \#(P\cap \t^2) = 1$, since $P$ and $F_1$ are
  connected normal surfaces. Thus, the choice of $P$ implies $P=F_1$,
  i.e., $P$ is fundamental with respect to $\t$ in $\p\setminus
  U(\t^0)$.\qed
\end{proof}

\subsection{Transformation into $\t_P$}
\label{sec:intoTP}

For any projective plane $P$ as in Lemma~\ref{lem:gibtP}, we define a
$2$--dimensional polyhedron
$$
Q_P = \big(\C^2\cap U(P)\big) \cup \d U(P).$$
Since $P$ is normal with respect to $\C$ and $\p\setminus U(P)$ is a
ball, $Q_P$ is the $2$--skeleton of a simple cellular decomposition of
$\p$, which we denote by $\C_P$. Any $2$--cell of $\C_P$ is contained
in the boundary of two different $3$--cells. Thus $\C_P$ is regular,
and the barycentric subdivision $\t_P=\C_P'$ of $\C_P$ is a
triangulation of $\p$.  The aim of this subsection is to relate $\t$
with $\t_P$ by a sequence of edge contractions and expansions.

We outline the construction of the transformation sequence.  All
ingredients are taken from~\cite{king1}--\cite{king3}, it is only
needed to adapt it to the present situation.
Since $B=\p\setminus U(P)$ is a ball, there is an embedding $H\co
S^2\times [0,1]\to B$ so that $H(S^2\times 0)=\d U(x)$ for some vertex
$x\in \C^0\subset \t^0$, and $H(S^2\times 1) = \d B = \d U(P)$.  
Let $c(H,\t^i)$ (resp.\ $c(H,\C^i)$) be the number of parameters
$\xi\in [0,1]$ for which the surface $H_\xi = H(S^2\times \xi)$ is not
in general position to $\t^i$ (resp.\ $\C^i$), for $i=1,2$. We assume
that $c(H,\t^1)$ is minimal.  An analysis of the Rubinstein--Thompson
algorithm as in~\cite{king1} yields an upper bound for $c(H,\C^1)$ in
terms of the number $n$ of tetrahedra of $\t$.
Techniques from~\cite{king2} allow to bound $c(H,\C^2)$ as well.  To
any surface $H_\xi$ that is in general position to $\C^2$, we define
an embedded $2$--complex $Q_\xi \subset \p$ that is the $2$--skeleton of a
simple regular cellular decomposition of $\p$, and $Q_1=Q_P$.
If $H_{\xi_0}$ is not transversal to $\C^2$ for some $\xi_0\in [0,1]$
and $\epsilon>0$ is sufficiently small, the complex
$Q_{\xi_0-\epsilon}$ is related to $Q_{\xi_0+\epsilon}$ by isotopy and
a bounded number of deletions and insertions of $2$--cells.
An application of Lemma~\ref{lem:contr} yields a transformation of
$\t$ into $\t_P$ by a bounded number of edge contractions and
expansions.

In the following lemma, \lq\lq normal\rq\rq\ shall mean \lq\lq
normal with respect to $\t$\rq\rq.
The estimate for $c(H,\t^1)$ is based on the construction of a
so-called \begriff{maximal normal sphere system}.  This is a system
$\Sigma\subset \p$ of disjoint normal $2$--spheres that are pairwise
not normally isotopic, so that any normal $2$--sphere in $\p\setminus
\Sigma$ is normally isotopic to a connected component of $\Sigma$.
\begin{lemma}\label{lem:Sigma}
  There is a maximal normal sphere system $\Sigma\subset \p\setminus U(P)$ 
  with at most $10n$ connected components and $\|\Sigma\| < 2^{181 n^2}$.
\end{lemma}
\begin{proof}
  We construct $\Sigma$ iteratively.  Define $\Sigma_1= \d U(P\cup
  \t^0)$.  Since $\d U(\t^0)$ meets each edge of $\t$ exactly twice and
  $\t$ has at most $2n$ edges, we have $\|\d U(\t^0)\| \le 4n$. 
  The projective plane $P$ is fundamental in $\p\setminus
  U(\t^0)$. Thus, by Theorem~\ref{thm:fundamental} and since $n\ge 24$, we
  have $$ \|\Sigma_1\| < 4n + 4n\cdot 2^{18 n} < 2^{19 n}.$$
  
  For $i\ge 1$, suppose that there is a connected component $N_i$ of
  $\p\setminus U(\Sigma_i)$ and a normal $2$--sphere $S\subset N_i$
  that is not normally isotopic to a connected component of
  $\Sigma_i$. It follows that $N_i$ is not a regular neighbourhood of
  $P$ or of a vertex of $\t$. We choose $S$ so that $\|S\|$ is
  minimal.
  
  Assume that $S$ can be represented as a sum $S_1+S_2$ of non-empty
  normal surfaces in $N_i$.  Since the Euler characteristic is
  additive and since there is no embedded projective plane in the
  $3$--ball $B=\p\setminus U(P)$, one of the summands, say $S_1$, is a
  sphere. It is not normally isotopic to a component of $\Sigma_i$,
  since otherwise $S_1+S_2$ would be the disjoint union of $S_1$ and
  $S_2$, thus would not be a sphere. We obtain a contradiction to the
  choice of $S$, since $\|S_1\|< \|S\|$. Thus $S$ is fundamental in
  $N_i$.
  
  We define $\Sigma_{i+1} = \Sigma_i \cup S$. By
  Theorem~\ref{thm:fundamental} and since $\|\d N_i\| \le \|\Sigma_i\|$,
  we have 
  $\|\Sigma_{i+1}\| \le \|\Sigma_i\| + \|\Sigma_i\|\cdot 2^{18n}$. 
  The iteration stops after at most $10n$ steps, by
  Theorem~\ref{thm:kneser}. Thus, we end with a maximal system $\Sigma$
  of normal $2$--spheres with
  $$
  \|\Sigma\| < \|\Sigma_1\| \cdot (2^{18n})^{10n-1} < 
                        2^{181 n^2}$$
  and at most $10n$ connected components.\qed
\end{proof}

\begin{lemma}\label{lem:TintoT_P}
  One can transform $\t$ into $\t_P$ by a sequence of less than $2^{184
    n^2}$ edge contractions and expansions.
\end{lemma}
\begin{proof}
  This lemma is a variant of Theorem~3 in~\cite{king3}.  We give here
  an outline of the proof, all details can be found in~\cite{king2} and~\cite{king3}.
  As in Lemma~34 of~\cite{king3}, there is an embedding $H\co
  S^2\times [0,1] \to \p\setminus U(P)$ in general position to $\t^1$
  so that $H(S^2\times 0)=\d U(x)$ for some $x\in \C^0\subset \t^0$,
  $H(S^2\times 1) = 2P$, and $c(H,\t^1)<
  \#(\Sigma)\cdot\|\Sigma\|\cdot 2^{18 n}$.
  Since $\C^1\subset \t^1$, by Lemma~\ref{lem:Sigma}, and since $n\ge
  24$, it follows
  \begin{eqnarray*}
    c(H,\C^1) &<& c(H,\t^1) \;< \;
    \#(\Sigma)\cdot\|\Sigma\|\cdot 2^{18 n}\\ &<& (10n)\cdot 2^{181 n^2 + 18
      n}
    \;<\;2^{182 n^2}.
  \end{eqnarray*}
  
  We denote $H_\xi=H(S^2\times\{\xi\})$ for $\xi\in [0,1]$.  By the
  choice of $P$ in Lemma~\ref{lem:gibtP}, both $H_1=2P$ and $H_0=\d
  U(x)$ are normal with respect to $\C$. In particular, for any
  $3$--cell $X$ of $\C$, any connected component of $H_0\cap \d X$
  (resp. of $H_1\cap \d X$) bounds a disc in $H_0\cap X$ (resp. in
  $H_1\cap X$).  Therefore $H$ satisfies the technical assumptions in
  Subsection~3.1 of~\cite{king2}.
  Hence by Lemmas~9 and~10 of~\cite{king2}, we can choose $H$ so that
  \begin{eqnarray*}
    c(H,\C^2) &\le& 1 + \chi\left(\C^2\cap U(P)\right)\underbrace{- \chi(\C^2) + \chi(\C^0)}_{< 10n}
                                            + c(H,\C^1)\\[-2mm]
        &<& 1 - \frac 12 \#(P\cap \C^1) + \hspace{.8cm}10\;n \hspace{.9cm}+ 2^{182 n^2}\\
        &<& 2^{183 n^2}.
  \end{eqnarray*}
  
  For $\xi\in [0,1]$, let $B^+(\xi)$ be the connected component of
  $\p\setminus H_\xi$ that contains $P$. We define $$Q_\xi = H_\xi\cup
  (\C^2\cap B^+(\xi)).$$
  If $H_\xi$ is in general position to $\C^2$
  then, by Lemma~13 in~\cite{king2}, $Q_\xi$ is the $2$--skeleton of a
  simple regular cellular decomposition of $\p$, whose barycentric
  subdivision is a triangulation $\t_\xi$.
  Let $\xi_0\in [0,1]$ so that $H_{\xi_0}$ is not in general position
  with respect to $\C^2$, and let $\epsilon>0$ be sufficiently small.
  It is shown after Lemma~13 in~\cite{king2} how one can transform
  $Q_{\xi_0-\epsilon}$ into $Q_{\xi_0+\epsilon}$. In the first step,
  depending on the type of non-transversality of $H_{\xi_0}$, one adds
  two $2$--cells to $Q_{\xi_0-\epsilon}$ with $2$ vertices in the
  boundary, or one adds one $2$--cell with at most $5$ vertices, or
  one deletes from $Q_{\xi_0-\epsilon}$ a $2$--cell with $3$ vertices.
  In the second step, one deletes a $2$--cell with at most $4$
  vertices.
  By an application of Lemma~\ref{lem:contr}, it follows that 
  $\t_{\xi_0-\epsilon}$ can be transformed into $\t_{\xi_0+\epsilon}$
  be a sequence of at most $22 + 18$ edge expansions or contractions.
  
  The complex $Q_0$ is obtained from $\C^2$ by insertion of a
  triangular $2$--cell. Hence one can transform $\t=\C'$ into $\t_0$
  by $14$ edge expansion. Furthermore, we have $\t_1=\t_P$.
  In conclusion, we obtain a sequence of less than
  $$
    14 + (22+18)\cdot c(H,\C^2) < 2^{184
    n^2}$$
  contractions and expansions relating $\t$ with $\t_P$.\qed
\end{proof}

\subsection{Standard triangulations}
\label{sec:standardtriang}

Let $P$, $\C_P$ and $\t_P$ be as in the preceding subsection.  Let
$\mathcal Z$ be a simple cellular decomposition of $P$ that is dual to
a triangulation.  We lift it along the fibres of the $I$--bundle
$\overline{U(P)}$ over $P$, and obtain a simple cellular decomposition
of $\overline{U(P)}$, so that each $d$--cell in $\d U(P)$ corresponds
to a $d$--cell of $\mathcal Z$, and the intersection of a $d$--cell in
$U(P)$ with $P$ is a $(d-1)$--cell of $\mathcal Z$, for $d=1,\dots,3$.
Since $\p\setminus U(P)$ is a ball, we obtain a simple cellular
decomposition of $\p$. It is easy to see that its barycentric
decomposition is a triangulation. We denote this triangulation by
$\t(\mathcal Z)$.

If there is  an edge $e$ in $\mathcal Z^1$ so that
$\mathcal Z^1 \setminus e$ is the $1$--skeleton of a cellular
decomposition of $P$ that is dual to a triangulation, then we replace 
$\mathcal Z$ by this simpler cellular decomposition.
\begin{figure}[htbp]
  \begin{center}
    {\psfrag{1}{\small$1$}\psfrag{2}{\small$2$}\psfrag{3}{\small$3$}
      \psfrag{4}{\small$4$}\psfrag{5}{\small$5$} \psfrag{Z}{$\mathcal
        Z_1$}\psfrag{Y}{$\mathcal Z_2$}\epsfig{file=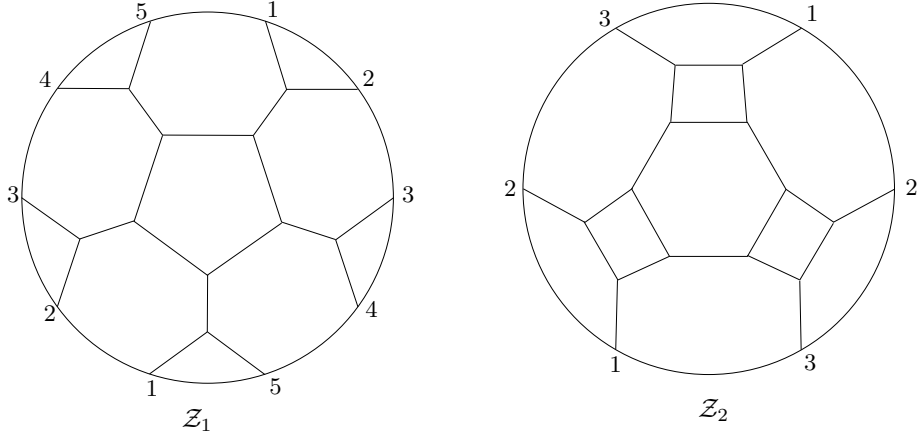,
        width=\linewidth}}
    \caption{The duals of the two irreducible simple triangulations of the 
        projective plane, according to Barnette}
    \label{fig:irredPtriang}
  \end{center}
\end{figure}
We iterate this process until it stops at a simple cellular
decomposition $\hat{\mathcal Z}$. It was proven by
Barnette~\cite{barnette} that $\hat{\mathcal Z}$ 
is one of the two decompositions $\mathcal Z_1,\mathcal Z_2$ depicted in
Figure~\ref{fig:irredPtriang} (opposite points in the boundary of the
discs are identified to obtain the projective plane $P$). 

\begin{lemma}\label{lem:standard}
  The triangulation $\t_P$ is related to one of the two standard
  triangulations $\t(\mathcal Z_1)$ and $\t(\mathcal Z_2)$ by a sequence
  of less than $2^{20 n}$ edge contractions.
\end{lemma}
\begin{proof}
  Let $\mathcal Z$ be the simple cellular decomposition of $P$ induced by
  $\C$. We have $\t(\mathcal Z) = \t_P$. 
  The deletion of edges of $\mathcal Z$ corresponds to the deletion of
  $2$--cells of $\C_P$. Each of these $2$--cells has four vertices in
  its boundary. Hence the deletion of one edge gives rise to $18$ edge
  contractions of $\t_P$, by Lemma~\ref{lem:contr}. 
  Since $\mathcal Z$ has $\frac 32\cdot\#(P\cap \C^1) \le \frac 32 \|P\| <
  2^{19n}$ edges, $\t_P$ is related to
  $\t(\hat{\mathcal Z})$ by a sequence of less then $18\cdot 2^{19n} <
  2^{20n}$ edge contractions.\qed
\end{proof}

\subsection{Proof of Theorem~\ref{thm:main}}
\label{sec:proofend}

Let $\t_1$ and $\t_2$ be two triangulations of $\p$ with at most $t$
tetrahedra. They are related to its barycentric subdivision $\t'_1$,
$\t'_2$ by a sequence of at most $5t$ edge expansions. Since $\t'_1$ and
$\t'_2$ have at most $12t$ tetrahedra, it follows from
Lemmas~\ref{lem:TintoT_P} and~\ref{lem:standard} that $\t'_1$ and
$\t'_2$ are related to one of the standard triangulations $\t(\mathcal
Z_1)$ and $\t(\mathcal Z_2)$ by sequences of less than
$$ 2^{184\cdot (12t)^2} + 2^{20\cdot 12 t} < 2^{26500 t^2}$$
edge expansions and contractions.

\begin{figure}[htbp]
  \begin{center}
    {\psfrag{1}{\small$1$}\psfrag{2}{\small$2$}\psfrag{3}{\small$3$}
    \psfrag{4}{\small$4$}\psfrag{5}{\small$5$} \psfrag{Y}{$\mathcal
      Z_1$}\psfrag{Z}{$\mathcal Z_2$}\epsfig{file=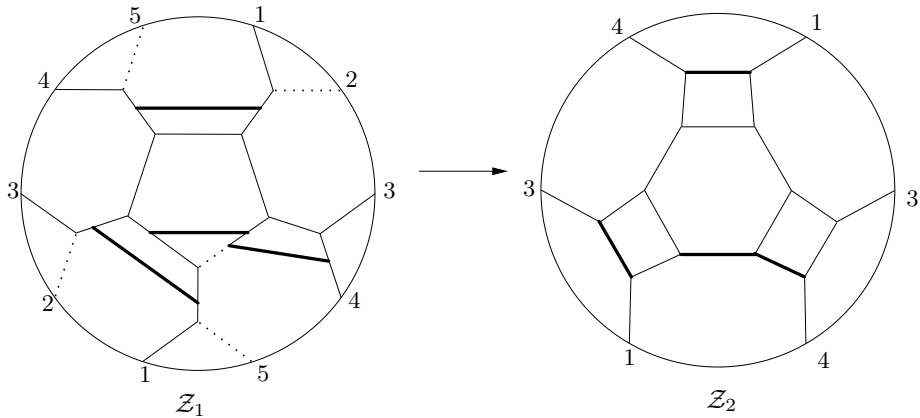,
      width=\linewidth}}
    \caption{Transforming $\mathcal Z_1$ into $\mathcal Z_2$}
    \label{fig:flips}
  \end{center}
\end{figure}
Figure~\ref{fig:flips} shows a transformation of $\mathcal Z_1$ into
$\mathcal Z_2$. One adds to $\mathcal Z_1$ four edges (the fat edges in the 
left part of Figure~\ref{fig:flips}) and deletes three edges of
the resulting cellular decomposition (the dotted edges in Figure~\ref{fig:flips}).
This corresponds to deletions and insertions of $2$--cells with four
vertices in the associated simple cellular 
decompositions of $\p$. Thus by Lemma~\ref{lem:contr}, $\t(\mathcal Z_1)$ is
related with $\t(\mathcal Z_2)$ by a sequence of $7\cdot 18 = 126$ 
edge expansions and contractions.

In conclusion, one can transform $\t_1$ into $\t_2$ by a sequence of
less than $126 + 10t + 2\cdot 2^{26500\cdot t^2} < 2^{27000\cdot t^2}$
edge expansions and contractions, which proves
Theorem~\ref{thm:main}.\qed

\end{document}